\documentclass[12pt, twoside, a4paper]{article}
\usepackage[english]{babel}  
\usepackage{float}
\usepackage[utf8]{inputenc}  
\usepackage[T2A]{fontenc}      
\usepackage{amsmath, amssymb, amsthm, amsfonts, amsthm, dsfont,mathrsfs,wrapfig,graphicx}
\usepackage[]{hyperref}
\usepackage{url}
\usepackage{epigraph}
\usepackage{accents}
\usepackage[LGR,T1]{fontenc} 

\DeclareSymbolFont{upgreek}{LGR}{cmr}{m}{n}
\SetSymbolFont{upgreek}{bold}{LGR}{cmr}{bx}{n}
\DeclareMathSymbol{\Eta}{\mathord}{upgreek}{`H}

\usepackage{abstract}
\usepackage{mathrsfs}
\usepackage[cal=cm,scr=dutchcal]{mathalfa}

\DeclareMathOperator{\supp}{supp}
\DeclareMathOperator{\dist}{dist}
\DeclareMathOperator{\pr}{pr}
\DeclareMathOperator{\const}{const}
\DeclareMathOperator{\clos}{clos}

\DeclareMathOperator{\diam}{{diam}}

\DeclareMathOperator{\card}{card}

\DeclareMathOperator{\spec}{spec}

\DeclareMathOperator{\mmod}{mod}

\DeclareMathOperator{\edge}{edge}
\DeclareMathOperator{\link}{link}
\DeclareMathOperator{\Cay}{Cay}
\DeclareMathOperator{\beg}{begin}
\DeclareMathOperator{\ennd}{end}
\DeclareMathOperator{\length}{length}

\DeclareMathOperator{\Ker}{Ker}
\DeclareMathOperator{\curl}{curl}
\DeclareMathOperator{\BS}{BS}

\newcommand{\eps}{\varepsilon}
\newcommand{\R}{{\mathbb R}}
\newcommand{\HH}{{\mathbb H}}

\newcommand{\NN}{\mathbb N}

\newtheorem{theorem}{Theorem}[section]
\newtheorem{lemma}[theorem]{Lemma}
\newtheorem{define}[theorem]{Definition}
\newtheorem{predl}[theorem]{Proposition}
\newtheorem{sled}[theorem]{Corollary}

\newcommand{\CC}{\mathbb C}

\usepackage{enumitem}

\usepackage{fancyhdr}
\renewenvironment{abstract}
{\small
	\begin{center}
		\bfseries \abstractname\vspace{-.5em}\vspace{0pt}
	\end{center}
	\list{}{%
		\setlength{\leftmargin}{2cm}
		\setlength{\rightmargin}{\leftmargin}%
	}%
	\item\relax}
{\endlist}

\makeatletter
\def\blindfootnote{\gdef\@thefnmark{}\@footnotetext}
\makeatother

\usepackage{mathtools}

\makeatletter
\def\blindfootnote{\gdef\@thefnmark{}\@footnotetext}
\makeatother
	
\begin{document}

\title{Coexact $1$-Laplacian spectral gap 
and exponential growth of a group}
\author{M. Dubashinskiy}
\date{\today}
\maketitle

	\blindfootnote{Chebyshev Laboratory, St.~Petersburg State University, 14th Line 29b, Vasilyevsky Island, Saint~Petersburg 199178, Russia.}
\blindfootnote{\hspace{0.2mm}e-mail: \href{mailto://mikhail.dubashinskiy@gmail.com}{\texttt{mikhail.dubashinskiy@gmail.com}}}
\blindfootnote{Research is supported by the Russian Science Foundation grant 19-71-30002.}

\blindfootnote{{Keywords:} Hodge Laplacian, exponential growth, spectral gap, $1$-cochain.
}

\blindfootnote{\hspace{0.0mm}{MSC 2020 Primary}: 58J50; Secondary:  20F65, 53C23.
}

\blindfootnote{This manuscript version is made available under the \href{http://creativecommons.org/licenses/bync-nd/4.0/}{CC-BY-NC-ND 4.0 license}.}

	
%

\thispagestyle{empty}

\renewcommand{\abstractname}{}	
\begin{abstract}
Let $\Gamma$ be a discrete finitely presented group. Pick any system $S$ of generators in $\Gamma$. In Cayley graph $\Cay(\Gamma)=\Cay(\Gamma, S)$ with edge set~$E$, glue with oriented polygons all the group relations translated to all the points of $\Gamma$; denote the obtained simply connected complex by $\Cay^{(2)}(\Gamma)$. 
We study non-negative \emph{Hodge--Laplace operator $\Delta_1$} on edge functions which is defined via complex $\Cay^{(2)}(\Gamma)$; $\Delta_1$ acts on
$$
\ell^2_{0,c}(E):= \clos_{\ell^2(E)}
\left\{\mbox{finitely supported closed $1$-(co)chains in }\Cay^{}(\Gamma)\right\}.
$$

We prove the following implication in the spirit of Kesten Theorem: \emph{if $\Delta_1|_{\ell_{0,c}^2(E)}$ has a spectral gap then $\Gamma$ either has exponential growth or is virtually $\mathbb Z$}. 
\end{abstract}

%
%

%

\bigskip

\bigskip

{\footnotesize
	\tableofcontents
}

\section{Introduction} 
\label{sec:Intro}

Let $G=(V,E)$ be countable oriented graph with degrees of vertices bounded from the above. Let $\tilde G$ be the non-oriented graph obtained from $G$ by forgetting the orientation of edges. 
Pick $D\in\mathbb N$ large enough. Consider all cycles in $\tilde G$ having lengths $\le D$. In $G$, glue each such cycle with a polygon. Choose any orientation of the latter polygons. 
We arrive to oriented $2$-dimensional complex, denote it by $G^{(2)}$ with implicit dependence on~$D$. Denote by $F$ the set of $2$-dimensional faces in $G^{(2)}$ which are polygons. Sometimes we write $F=FG$ and also $E=EG$ to indicate the dependence of these sets on $G$. Any of sets $V, E, F$ is endowed with counting measure which we denote by $\card$. In graph $G$, we define graph metric $\dist_G$ at $V\cup E$ along edges in $E$ so that any edge has length $1$.

If $\Gamma$ is a finitely generated group, $S$ is any of its generating sets (symmetrized or not) then we may consider $G=\Cay(\Gamma,S)$, Cayley graph of $\Gamma$; then $V=\Gamma$. If $\Gamma$ is also finitely presented, that is, given by a finite number of relations then we assume that $D$ in the definition of $G^{(2)}=\Cay^{(2)}(\Gamma,S)$ is $\ge$ than length of any of the defining relations. For general $G$, we assume that $D$ is such that 
\begin{equation*}\mbox{$G^{(2)}$ is simply connected}
\end{equation*}
(and that such $D$ does exist).


A 
\emph{$k$-cochain}, $k$ is $0$, or $1$, or $2$, is a function from $V$, or $E$, or $F$, respectively, to $\R$. We often understand  cochains as chains.
Discrete differentials (coboundaries)
$$
\{0\mbox{-cochains}\}\xrightarrow{d}\{1\mbox{-cochains}\}\xrightarrow{d}\{2\mbox{-cochains}\}
$$
and boundary operators 
$$
\{2\mbox{-cochains}\}\xrightarrow{\partial}\{1\mbox{-cochains}\}\xrightarrow{\partial}\{0\mbox{-cochains}\}
$$
are introduced in the standard way with respect to the orientation of edges and faces.  Since valencies of vertices are bounded, all these operators are also bounded with respect to $\ell^2$-norms on cochains. 
We have $(d|_{\ell^2(V)})^*=\partial|_{\ell^2(E)}$, $(d|_{\ell^2(E)})^*=\partial|_{\ell^2(F)}$. Indeed, discrete integration by parts is valid for finitely supported cochains and is proved for  $\ell^2$-cochains by $\ell^2$-approximation. 

If $\gamma$ is an oriented path in $G$ then we may define $1$-(co)chain $f_\gamma$: for $e\in E$, let $f_\gamma(e)$ be the number of passes of $\gamma$ through $e$ in its direction minus number of passes of $\gamma$ over $e$ in its reversed direction. 
Then we have $\partial f_\gamma=0$.

Our space of interest is
\begin{equation*}
	\ell_{0,c}^2(E) := \clos_{\ell^2(E)}\{f\colon E\to \R\mid \partial f=0, \, \supp f
	\,\mbox{is finite}\}.
\end{equation*}
Any of $1$-cochains at the right-hand side can be (convexly) decomposed into simple finite loops. Thus, $\ell_{0,c}^2(E)$ is $\ell^2$-closed linear span of (co)chains of the form $f_\gamma$ with $\gamma$ a finite loop in $G$.

Laplace operator on $0$-cochains is 
$$
-\Delta_0=\partial d\colon \left(\mbox{functions on }V\right)\to \left(\mbox{functions on }V\right).
$$
A discrete integration by parts leads to the following Hodge-type decomposition: 

\begin{predl}
	
	We have 
	\begin{equation}
		\label{eq:orth_decomp}
		\{f\in\ell^2(E)\colon \partial f=0\}= \ell_{0,c}^2(E)\oplus_{\ell^2(E)} \{du\mid u\colon G\to \R, \, \Delta_0 u = 0, \, du\in \ell^2(E)\}.
	\end{equation}
\end{predl}

The second summand in the right-hand side of the latter relation is \emph{$\ell^2$-cohomology} of $G$. It is known to be invariant with respect to change of generating system in a group: 
the factorspace nature of cohomology allows to implement "discrete change of variables" from one to another set of generators. Cohomology is invariant with respect to more general quasiisometries. 

Now, we pass to spectral estimates for $1$-cochains. Define non-negative Laplacian operator $\Delta_1:=\partial d+d\partial \colon \ell^2(E)\to \ell^2(E)$. On $\ell^2_{0,c}(E)$, our space of interest, this reduces to $\partial d$. 

We have one more Hodge-type decomposition:
$$
\ell^2(E) = \clos_{\ell^2(E)}\{du\mid u \colon V\to \R, \, \supp u \mbox{ is finite}\} \oplus_{\ell^2(E)} \{f\in\ell^2(E)\colon \partial f=0\}.
$$
Spectral questions for $\Delta_1$ on the first summand are generally reduced to the same for $\Delta_0$ on $\ell^2(V)$.   What concerns decomposition (\ref{eq:orth_decomp}) for $\{f\in\ell^2(E)\colon \partial f=0\}$, operator $\Delta_1$ vanishes at the second summand of its right-hand side, $\ell^2$-cohomology. Also, by the definition of $\ell^2_{0,c}(E)$ and by $\ell^2$-approximation, we see that $\Delta_1(\ell^2_{0,c}(E))\subset\ell^2_{0,c}(E)$.

\begin{define}
	We say that $\Delta_1$ \emph{has a spectral gap at $\ell^2_{0,c}(E)$} \emph{(}or just that graph $G$ has \emph{{coexact $1$-Laplacian spectral gap}}\emph{)} if $$\spec\left(\Delta_1|_{\ell^2_{0,c}(E)}\right)\cap [0,\eps)=\varnothing$$ for some $\eps>0$ small enough.
\end{define}

Applying discrete integration by parts, we conclude that this is equivalent to the estimate 
\begin{equation}
	\label{eq:rot_estim}
\langle f, f \rangle_{\ell^2(E)}	\le 1/\eps\cdot \langle df, df\rangle_{\ell^2(F)}	
\end{equation}
for $f\in \ell^2_{0,c}(E)$. It is enough to check the latter 
only for finitely supported closed $1$-cochains $f$. Also, we conclude that if $1$-Laplacian has a coexact spectral gap then it will be so if we enlarge $D$ in the construction of $G^{(2)}$ or glue some extra faces to $G^{(2)}$ in a locally finite manner. 
 
 We state quasiinvariance result as below, with possibility to add only edges. It seems feasible to preserve spectral gap under more general quasiisometric transformations of a graph, the ones with possibility to add or remove vertices in a locally finite way. 
 
 \begin{predl}
 	\label{predl:general_invariance} 
 	Let $G_1=(V,E_1)$ and $G_2=(V,E_2)$ be graphs with the same vertex set $V$ and with $E_1\subset E_2$. 
 	Assume the following:
 	\begin{enumerate}
 		\item $G_1$ is connected;
 		
 		\item degrees of vertices in $G_2$ are bounded from the above;
 		
 		\item $\sup\limits_{e\in E_2\setminus E_1} \dist_{G_1}(\beg e, \ennd e) < +\infty$ with the obvious notation.  
 		
 	\end{enumerate}
	\emph{(}In other words, metrics $\dist_{G_1}$ and $\dist_{G_2}$ on $V$ are bilipshitz equivalent.\emph{)} Then, $G_1$ has a coexact $1$-Laplacian spectral gap  \emph{(}with some $D$ implied at the construction of $G_1^{(2)}$\emph{)}
	if and only if  $G_2$ has such a spectral gap \emph{(}with some $D'$ for $G_2^{(2)}$\emph{)}.
 \end{predl}
 
 \noindent Proof of this quasiinvariance  is similar to the proof of quasiiinvariance of $\ell^2$-cohomology; both are based on orthogonal projection. We give a detailed argument in the Appendix. Notice also that the proof is constructive: we may estimate $D'$ via $D$ and the supremum from the third assumption of Proposition \ref{predl:general_invariance}, and vice versa.
 
 
 \begin{sled}
 	For two Cayley graphs of the same finitely-presented group but with different generating sets, properties of existence of a coexact $1$-Laplacian spectral gap on them  are equivalent.
 \end{sled}

Our main result is

\begin{theorem}
	\label{th:exp_growth}
	Let $\Gamma$ be a countable finitely presented group.
	If $\Delta_1$ has a spectral gap on 
	$\ell_{0,c}^2(E)$ 
	then either $\Gamma$ has exponential growth, 
	 or $\Gamma$ is virtually infinite cyclic. 
\end{theorem}

Let us recall the well-known Kesten Theorem on Laplacian on vertices of a graph:

\begin{theorem}[\cite{K59}]
	\label{th:Kesten}
	Let $\Gamma$ be a finitely generated group. Then, for $0$-Laplacian in $\Cay(\Gamma)$,  $$0\in \spec\left(-\Delta_0|_{\ell^2(\Gamma)}\right)$$ if and only if $\Gamma$ is amenable.
\end{theorem}

%

Non-amenability of a group, that is, the existence of a spectral gap for $-\Delta_0$, easily implies exponential growth. The reverse is not true, in general. Thus, it is natural to ask, for example, whether Baumslag--Solitar groups $\langle\mathscr a, \mathscr b\mid \mathscr b^{-1}\mathscr a\mathscr b=\mathscr a^n\rangle$, $n\in\mathbb N$, have a spectral gap for $1$-Laplacian. These groups are \emph{non-elementary amenable} but have exponential growth. Such groups are not covered by Theorem \ref{th:exp_growth}, and it is still unclear for the author whether $1$-Laplacian has a spectral gap on them.

If we assume the contrary to Theorem \ref{th:exp_growth}, then, first, $\Gamma$ cannot have two \emph{ends} since in this case $\Gamma$ is virtually cyclic, see, e.g., \cite{Me08}. Second, $\Gamma$ cannot have infinitely many ends because then $\Gamma$ has exponential growth; the latter follows from Stallings Theorem and from results of \cite{HB00} but, of course, can be proved directly. So, by Freudenthal--Hopf Theorem we may assume that $\Gamma$ has one end.

For $\mathscr L>0$, denote by $\mathbb T_{\mathscr{L}}$ a circle of length $\mathscr L$. On $\mathbb T_{\mathscr L}$, one may measure distances along this loop. To prove Theorem \ref{th:exp_growth}, we need the following

\begin{lemma}[on loop embedding]
\label{lemma:embed}
Suppose that $\Gamma$ has subexponential growth and just one end.

Let  $C>20$, $x\ge 1$. 
There exist 
$\mathscr L>2DCx$ and injective naturally parametrized $\gamma\colon \mathbb T_{\mathscr L}\to \Cay(\Gamma)$ such that, for $t_1, t_2\in\mathbb T_{\mathscr L}$ with $\gamma(t_1), \gamma(t_2)\in \Gamma$, 
\emph{
\begin{equation}
	\label{eq:x_bilip}
\mbox{if }\dist_{\Cay(\Gamma)}(\gamma(t_1), \gamma(t_2)) \le x 
\mbox{ then }
\dist_{\mathbb T_{\mathscr L}}(t_1, t_2) \le Cx.
\end{equation}
}
\end{lemma}

In fact, we are able to make $\mathscr L$ arbitrarily large with fixed $x$.

Now, let us briefly recall the proof of Kesten Theorem \ref{th:Kesten} to compare it to our argument. Non-amenability of $\Gamma$ means that
\begin{equation}
\label{eq:isoperimetr}
\|du\|_{\ell^1(E)} \gtrsim \|u\|_{\ell^1(\Gamma)}
\end{equation}
for $u=\mathds 1_E$, $E$ ranges all finite subsets in $\Gamma$. By discrete version of coarea formula, this is equivalent to the same for any finitely supported $u\colon \Gamma\to\R$. 
The spectral gap condition $0\notin \spec\left(-\Delta_0|_{\ell^2(\Gamma)}\right)$ means that  
$$
\langle u, u\rangle_{\ell^2(\Gamma)} \lesssim \langle du, du\rangle_{\ell^2(E)}, ~~~ u\colon\Gamma\to\R \mbox{ is finitely supported}.
$$
To obtain this from (\ref{eq:isoperimetr}), it remains to insert $u^2$ instead of $u$ to (\ref{eq:isoperimetr}) and apply Cauchy--Bunyakovsky--Schwartz inequality.

In the first step of the latter argument, we assemble a function $u\colon \Gamma\to\R$, say, non-negative one, from its super-level sets $\mathds 1_{\{u\ge \mathscr t\}}$, $\mathscr t$ ranges $[0,+\infty)$; we also assemble $du$ from $d\mathds 1_{\{u\ge \mathscr t\}}$. (Both decompositions are $\ell^1$-convex.) Thus, in Kesten Theorem, we deal with "sets of codimensions $0$ and $1$". At least, we will have such genuine codimensions in the case of a manifold instead of a group, the corresponding result linking spectra and isoperimetry is known as Cheeger--Yau inequality, see \cite{Cheeger}, \cite{Yau}.

Unlike this, in our argument we work with dimension $1$ sets --- loops, in particular, as in Lemma \ref{lemma:embed}. Also, in Section \ref{sec:spectra_homology} we bound $1$-cycles with $2$-dimensional surfaces. 

Notice also that an analogue of Cheeger--Yau inequality for $1$-forms was obtained in \cite{BC22} in the case of manifolds. 
Coexact $1$-Laplacian spectrum is indeed related to appropriate isoperimetric ratio, namely, to $\sup_{\gamma}\inf_h|h|/\length\gamma$ with $\gamma$ ranging homologicaly trivial loops at a manifold $\mathscr M$ and $h$ be a $2$-dimensional chain in $\mathscr M$ bounding $\gamma$; here, $|h|$ is area of $h$. Some Poincar\'e-type estimates for operator $d$ on coclosed $1$-forms are  possible if isoperimetric ratios as above are bounded from the below. But, in \cite{BC22}, authors impose the condition of finite diameter of $\mathscr M$  which is not our case; also, \cite{BC22} does not deal with effects of negative curvature.

What concerns spaces with negative curvature, let us mention recent works \cite{AAGLZ24}, \cite{R23} devoted to $3$-dimensional hyperbolic manifolds. It turns out that, first, $1$-coexact spectral gap is related to exponential growth of torsion $1$-homology of the manifolds; second, there are relations  between the spectral gap and isoperimetric ratios. The latter estimates from \cite{R23} are also volume-dependent, as in \cite{BC22}.

\medskip

This paper is organized as follows. In Section \ref{sec:maueutic}, we explain our interest to 
the study of $1$-Laplacian spectra. This Section is not used in the proof of Theorem \ref{th:exp_growth}. 
In Section~\ref{sec:loop}, we prove Lemma \ref{lemma:embed} by dropping lots of geodesic perpendiculars in a branching way. In Section \ref{sec:spectra_homology}, we conclude the proof of Theorem~\ref{th:exp_growth}. This is done by approximating the resolvent $\Delta^{-1}_1$ by polynomials of $\Delta_1$ provided that $0\notin\spec\Delta_1|_{\ell_{0,c}^2(E)}$. Next, we apply this to $1$-cochain given by curve $\gamma$ from Lemma \ref{lemma:embed}. We put some metric control on the approximating process as implemented in $\Cay(\Gamma)$ and also  make use of homological nature of $\Delta_1$: this operator is divisible by $\partial\colon \ell^2(F)\to\ell^2(E)$ at $\ell_{0,c}^2(E)$.  Finally, in Section \ref{sec:examples} we check the most natural examples of Cayley graphs. 

\medskip

\noindent {\bf Some notation.} For a set $A$ we denote by $\card A$ the number of its elements. 
If $v_1, v_2$ are vertices of some oriented graph then we denote by $\edge(v_1, v_2)$ the oriented edge in the graph under consideration provided that the edge exists. If $e$ is an edge in a oriented graph or $\gamma$ is an oriented path in a metric space then we write $\beg e$, $\beg\gamma$ for their beginnings and $\ennd e$, $\ennd\gamma$ for their endpoints, respectively. The notation $\length\gamma$ is obvious. 

We write $\mathcal B_X(x, \rho)$ for the open ball in a metric space $X$ centered in a point $x\in X$ and having radius $\rho\ge 0$.

\section{Motivation: functional-analytic approach}
\label{sec:maueutic}

Let us expose some considerations lead author to the study of $\spec\Delta_1$ on $\ell^2_{0,c}(E)$. Reader may skip this Section safely until Section \ref{sec:examples}. 

The space $H=\ell_{0,c}^2(E)$ is usually of infinite countable dimension. 
All such Hilbert spaces $H$ are isomorphic, and, abstractly speaking, there is nothing to classify. Instead, we may try to classify tuples $(H, f_1, f_2, \dots)$ with $\{f_1, f_2, \dots\}$ is a countable system in an abstract Hilbert space $H$. We may ask for a classification up to action of $GL(H)$; the latter is the group of all linear bounded, boundedly invertible operators in $H$. 

If $f_j$ as above are of geometric nature then we may impose geometric restriction on them. For example, if $H=\ell_{0,c}^2(E)$  then we may require $\sup\limits_{j\in\mathbb N} \diam\supp f_j<+\infty$.

Also, if some group $\Gamma$ acts on $H$ then we may require that set $\{f_j\}_{j\in\mathbb{N}}$ is $\Gamma$-invariant, or up to finite index subgroup, or consists of finite number of orbits; or even that it is an infinite union of orbits with limited growth of supports, etc.

We return to functional-analytic restrictions on $\{f_j\}_{j\in\mathbb N}$. What concerns orthonormal systems in (co)homology, author would be amazed by an example of such a basis which is not an eigenbasis of a self-adjoint operator nor is obtained by Gram--Schmidt process.

We may relax orthonormality condition. 
Recall that an image  of an orthonormal basis under an action of an operator from $GL(H)$ is called a  \emph{Riesz basis} in $H$. So, we may ask for an existence, say, of a localized equivariant Riesz bases in $\ell^2_{0,c}(E)$. From the first glance, action of $GL(H)$ seems to be an adequate functional-analytic counterpart of procedure of change of a generating system in a group since $GL(H)$ is "softer" than the group of unitary operators on $H$.

Alas, Riesz basis condition still seems to author to be too rigid in our topological setting: generally, we have just rare topological spaces with clear \emph{basis} even in the usual unnormed homology space (with coefficients in $\R$). We meet such example in Section \ref{sec:examples} (standard hyperbolic plane tilings), 
see also \cite{D16} for planar disk with holes. What concerns the standard procedure of retracting a graph onto a bouquet of circles, it generally does not lead to well-localized $\R$-homology bases, as we wished before; nor does it automatically lead to group-equivariant bases in the case of presence of a group action.

Also, one should immediately raise the question on invariance of existence of good Riesz bases with respect to, say, change of generators in group. If $\Gamma$ is a group, $S$ is any of its generating sets then, having an equivariant well-localized Riesz basis for $\Cay(\Gamma, S)$ we easily construct such a basis for $\Cay(\Gamma, S\cup\{s\})$ for any $s\in\Gamma$. But it is completely unclear for the author how reconstruct Riesz bases under  removal of a generator. 


So, a property to be a basis in $1$-homology is too rigid, even without Riesz condition. Instead 
we consider the notion of a \emph{frame} which turns to be more flexible.

\begin{define}
	Let $H$ be a separable Hilbert space \emph(over $\R$ or $\CC$\emph ), $f_1, f_2, \dots \in H$. We say that $\{f_n\}_{n\in\NN}$ is a \emph{frame} in $H$ if it satisfies the following \emph{almost Parseval} condition: for any $g\in H$, we have
	\begin{equation*}
		C^{-1} \|g\|^2_H \le \sum\limits_{n\in\mathbb N} |\langle g, f_n\rangle|^2 \le C \|g\|^2_H
	\end{equation*}
	with some $C\in (0,+\infty)$ not depending on $g$.
\end{define}


From (\ref{eq:rot_estim}) we conclude that the following assertions on a group $\Gamma$ are equivalent:
\begin{itemize}
	\item $\Delta_1$ has a spectral gap on $\ell^2_{0,c}(E)$;
	\item when $\gamma$ ranges the set of all oriented loops in $\Cay(\Gamma)$ of length $\le D$, family $\{f_\gamma\}$ of $1$-(co)chains generated by loops is a frame in $\ell^2_{0,c}(E)$.
\end{itemize}

Indeed, upper estimate from frame definition is immediate provided that degrees of vertices are bounded from the above.

\begin{predl}
\label{predl:scf_def_equivalence}
	Let $G$ be a graph with degrees of vertices bounded from the above. Suppose that there exists some frame $\{g_n\}_{n\in\NN}$ in $\ell^2_{0,c}(E)$ with $$D:=\sup\limits_{n\in\NN}\diam\supp g_n < +\infty.$$ Then $\Delta_1$ has a spectral gap in $\ell^2_{0,c}(E)$ with the same $D$ implied in the construction of~$G^{(2)}$.
\end{predl}

Proof 
is elementary and is given at the Appendix.

To conclude this Section, we just mention that functional-analytic viewpoint is applicable also to linking properties of two lattices $\mathbb Z^3$ and 
$\mathbb Z^3+(1/2,1/2,1/2)$ in $\R^3$; denote the corresponding graphs by $(V_1, E_1)$ and $(V_2, E_2)$. For any two finitely supported cycles $\gamma_1$ in $(V_1, E_1)$ and $\gamma_2$ in $(V_2, E_2)$, one defines \emph{linking number} $\link(\gamma_1, \gamma_2)$ which measures how much times $\gamma_2$ wires around $\gamma_1$. Thus, any finitely supported $f\colon E_1\to \R$ with $\partial_{(V_1,E_1)}f=0$ gives a functional $\link(f, \cdot)\colon \ell^2_{0,c}(E_2)\to \R$. One may study functional-analytic properties of such functionals when $f$ is well-localized. We do not proceed this here.

%
%


\section{Loop embedding} 
\label{sec:loop}

Here, we prove Lemma \ref{lemma:embed}. Fix  $x$ and assume that conclusion of Lemma \ref{lemma:embed} is not valid for this $x$.

If $a, b\in\Gamma$ then we denote by $[a,b]$ a geodesic segment joining $a$ and $b$; if there is a plenty of such shortest paths then $[a,b]$ can be either specified explicitly or is taken in an arbitrary way. 
We will not arrive to an ambiguity. 
If also $c\in\Gamma$ then denote by $[a,b]+[b,c]$ concatenation of two such segments passed from $a$ to $c$.

For $a, b, c_1, c_2\in \Gamma$ and a geodesic segment $[c_1, c_2]$ such that $b\in[c_1, c_2]$, we say that $[a,b]$ is a \emph{perpendicular} to $[c_1, c_2]$ if $[a,b]$ is the shortest path from $a$ to a point at $[c_1, c_2]$ or one of such paths if there is lots of them; in this case, we write $[a,b]\bot[c_1, c_2]$. In Lemmas below we often have $b=c_1$ and thus may speak about perpendicular angles.

The following Lemma shows that we may use something besides geodesic segments to satisfy (\ref{eq:x_bilip}) locally, namely, that perpendicular angles are also useful to this end (though they are not closed).

\begin{lemma}
\label{lemma:perpendicular} Let $a, b, c\in\Gamma$, $[a,b]$ and $[b,c]$ be any of geodesic segments with given ends. Suppose that  $[a,b]\bot [b,c]$.

If $\length([a,b]+[b,c]) > Cx$ then $\dist_{\Cay(\Gamma)}(a,c)>x$. 
\end{lemma}

\noindent {\bf Proof.}
Indeed, if $\dist_{\Cay(\Gamma)}(a,c)\le x$ then $\length[a,b]\le x$ since $[a,b]\bot[b,c$]. By triangle inequality, $\length[b,c]\le 2x$. Therefore, $\length([a,b]+[b,c]) \le Cx$ if $C>3$.
$\blacksquare$

\medskip

Assume that $\Gamma$ does not have exponential growth. Pick $N\in\mathbb N$ such that $$\card\mathcal B_{\Gamma}(1_\Gamma, 4CDxN) < 2^{N-1}.$$

We construct an almost-binary tree ordered by levels with vertex set $\mathcal T$ as follows. Its root $\nu_0\in\mathcal T$ has level $0$ and it has one \emph{right} descendant which we denote by $\nu_0r$, the latter is at level $1$. Any other vertex $\nu\in\mathcal T$ at level $\le N-1$ besides the root has one left and one right descendant at the next level, denote them by $\nu l$ and $\nu r$, respectively;  vertices at level $N$ are leafs and do not have descendants. Our notation allows us to write $\nu w\in\mathcal T$ for $\nu\in\mathcal T$ and $w$ a word in alphabet $\{l,r\}$ which is not too long.

We are going to construct a mapping $\phi \colon \mathcal T\to \Gamma$, and also, 
for any $\nu\in\mathcal T$ of level $\le N-1$, an oriented geodesic segment $\lambda(\nu)$ with the following properties.

\begin{itemize}
	\item  For adjacent $\nu_1, \nu_2\in\mathcal T$, we have:
\begin{equation}
	\label{eq:vary_step}
	\dist_{\Cay(\Gamma)}(\phi(\nu_1), \phi(\nu_2))\in [2CDx, 4CDx].
\end{equation}
	\item Let $\nu\in\mathcal T$ be a vertex at level $k\in[0,N-1]$. Then segment $\lambda(\nu)$ starts at $\phi(\nu)$ and then passes $\phi(\nu r), \phi(\nu rl),\phi(vrl^2), \dots, \phi(\nu rl^{N-k-1})$ in the order from $\phi(\nu)$ to $\phi(\nu rl^{N-k-1})$ and may also pass some other vertices from $\Gamma$ between them.
	
	\item Let $\nu\in\mathcal T$ be not a leaf and let  $\nu'$ be any of $\nu r, \nu rl, \nu rl^2, \dots, \nu rl^{N-k-2}$. Then $\lambda(\nu')$ reversed is perpendicular to $\lambda(\nu)$ (at the point $\phi(\nu')\in\lambda(v)\cap\lambda(v')$, by our construction).
	
	\item For technical reasons, we ask that segments $\lambda(\nu)$ can be taken arbitrarily long. We will specify the choice of their lengths below in a back-recursive way.
\end{itemize}

This construction is rather clear in terms of segments $\lambda(\nu)$. Assume, for a while, that we are extremely lucky and are able to build long enough perpendiculars to any geodesic segment at any of its point. Then we construct segments $\lambda(\nu)$ in the following order. For $\lambda(\nu_0)$ we take a long enough geodesic segment, and take its starting point for $\phi(\nu_0)$. Assume that, for $\nu\in\mathcal T$ of level $k\in[0,N-1]$, segment $\lambda(\nu)$ long enough is already built. Pick $\phi(\nu r), \phi(\nu rl), \dots, \phi(\nu rl^{N-k-1})$ lying at $\lambda(\nu)$ in this order to have
\begin{multline}
\label{eq:fixed_step}
\dist_{\Cay(\Gamma)}(\phi(\nu), \phi(\nu r))=
\dist_{\Cay(\Gamma)}(\phi(\nu r), \phi(\nu rl))=\dist_{\Cay(\Gamma)}(\phi(\nu rl), \phi(\nu rl^2))=\\=\dots=\dist_{\Cay(\Gamma)}(\phi(\nu rl^{N-k-2}), \phi(\nu rl^{N-k-1}))=\lceil2CDx\rceil.
\end{multline}
Further, for $j=0,1,  \dots, N-k-2$, let $\lambda(\nu rl^j)$ be long enough  perpendicular to $\lambda(\nu)$ built at point $\phi(\nu rl^j)$. Then repeat our procedure for newly constructed segments and stop at leafs of level $N$ in $\mathcal T$.

By the choice of $N$, there are two leafs $\nu_1, \nu_2\in\mathcal{T}$ with $\phi(\nu_1)=\phi(\nu_2)$, $\nu_1\neq \nu_2$. This does not immediately lead to construction of curve $\gamma$ from Lemma \ref{lemma:embed}. We thus also need to apply a loop shrinking procedure as in Lemma \ref{lemma:shrinking}.

But, generally, we are not able to build perpendiculars to any geodesic at any prescribed points. We return to formal consideration and start with addressing the questions on perpendiculars: either they do exist up to shifting the basepoint by the distance $\le 2CDx$, or $\Gamma$ has two ends. Otherwise, we construct curve $\gamma$ for Lemma \ref{lemma:embed} if some  
natural auxiliary steps fail for $\Gamma$.


\begin{lemma}
	\label{lemma:perpendicular_stab}
	Let $a_1, a_2, b, c\in\Gamma$, $a_1$ and $a_2$ adjacent in $\Cay(\Gamma)$. Suppose that $[a_1, b], [a_2,c]\bot[b,c]$. 
	
	If $\length[b,c]>2DCx$ then we may construct a loop from Lemma \ref{lemma:embed} for our $x$.
\end{lemma}

\noindent {\bf Proof.}
Consider all pairs of points $(a_1', a_2')$ with $a_1'\in[a_1,b]$, $a_2'\in[a_2,c]$ such that $\dist_{\Cay(\Gamma)}(a_1', a_2')\le x$. Among all such pairs, take the one "closest" to $[b,c]$, namely, the one with minimal $\length[a_1',b]+\length[a_2',c]$. Let  $\gamma$ be natural parametrization of loop $[a_1', b]+[b,c]+[c,a_2']+[a_2', a_1']$. We claim that this $\gamma$ satisfies all the conditions from Lemma \ref{lemma:embed}. Its length $\mathscr L$ is $>\length[b,c]>2CDx$.

By the choice of $a_1', a_2'$, for any two consequent sides of the geodesic quadrilateral~$\gamma$, one of them is perpendicular to another.  
Thus, if two points belong to adjacent sides of~$\gamma$ then (\ref{eq:x_bilip}) for such points follows from Lemma \ref{lemma:perpendicular}. Injectivity for such pairs of points also is immediate.

By the choice of $a_1'$, $a_2'$, we have $\dist_{\Cay(\Gamma)}([a_1',b], [a_2',c])=x$,  with strict distance minimum attained at $a_1'$ and $a_2'$. This implies (\ref{eq:x_bilip}) for $\gamma(t_1)\in [a_1',b]$, $\gamma(t_2)\in [a_2',c]$ in the notation of Lemma \ref{lemma:embed}; also injectivity for such pairs follows.

Finally, assume that $\dist_{\Cay(\Gamma)}([a_1', a_2'],[b,c])\le x$. Then $$\dist_{\Cay(\Gamma)}(a_1', [b,c]), \dist_{\Cay(\Gamma)}(a_2', [b,c]) \le 2x.$$ Since $[a_1',b], [a_2',c] 
\bot [b,c]$, we derive that $\length[a_1',b], \length[a_2',c]\le 2x$. Therefore, 
$$
5x \ge \length([b, a_1']+[a_1', a_2']+[a_2',c])\ge \length[b,c]\ge 2CDx$$
which is impossible for $C>5$.
$\blacksquare$

\medskip

In the construction of the almost-binary tree mapping, we need the following 

\begin{lemma}
	\label{lemma:perpendicular_stab2}
	Let $a,b,c\in\Gamma$ such that $[a,b]\bot[b,c]$. 
	Let $a_1$ be point on $[a,b]$ closest to $c$.
	
	If $\length[a_1,b]>2CDx$ then it is possible to construct a curve required in Lemma~\ref{lemma:embed} for our $x$.
\end{lemma}

\noindent {\bf Proof.} Similar to the proof of Lemma \ref{lemma:perpendicular_stab}, $C>6$ is enough.
$\blacksquare$

\medskip

Now let us prove the possibility to build geodesic perpendicular of length $L\ge 1$ not far enough from a given point at a geodesic segment.

\begin{lemma}
	\label{lemma:build_perpendicular}
	Let $C, D, x, N$ be fixed.  For any $L\ge 1$ there exists $L'\ge 10CDxN$ large enough with the following property.
	
	Let $\lambda$ be a geodesic segment  in $\Cay(\Gamma)$ starting at a vertex $v_0\in\Gamma$. Assume that:
	\begin{itemize}
		\item $\length\lambda \ge  L'$;
		\item there exists a geodesic segment $\lambda_1$ ending at $v_0$ such that 
		$\length\lambda_1\ge L'$ and $\lambda$ reversed is a perpendicular to $\lambda_1$. 
	\end{itemize}
	Denote by  $\lambda_2$ the subsegment of $\lambda$ starting at $v_0$ and of length $10CDxN$.

	Then, for any subsegment $\lambda_3$ in $\lambda_2$ with $\length\lambda_3=2CDx$, there exists a segment $[a,b]$ with $a\in\Gamma$, $b\in\lambda_3$, $\length[a,b] >L$ and such that $[a,b]\bot \lambda_2$. Otherwise, either $\Gamma$ has $\ge 2$ ends, or we success in constructing curve $\gamma$ for $x$.
	
	Length $L'$ depends not only on $C, D, x, N, L$ but also on $\Gamma$ if it is a 
	group with one end.
\end{lemma}

\noindent {\bf Remark.} Assume $\lambda_0$ is a geodesic segment in $\Cay(\Gamma)$ of even length $2L'$ with $\beg\lambda, \ennd\lambda\in\Gamma$. If $\lambda$ is any of halves of $\lambda_0$ starting in its  middle point then the assumption of Lemma \ref{lemma:build_perpendicular} is valid for this $\lambda$. This is because, in our terms, flat angle is also a right angle, and we may take the rest half of $\lambda_0$ for $\lambda_1$.

\medskip

\noindent {\bf Proof of Lemma \ref{lemma:build_perpendicular}.} First, let $\lambda_2$ range the family of all the geodesic segments in $\Cay(\Gamma)$ of length $10CDxN$. Let $U_{\lambda_2}$ be $L$-neighborhood of $\lambda_2$ in $\Cay(\Gamma)$. Consider sets $\Cay(\Gamma)\setminus U_{\lambda_2}$. If $\Gamma$ has just one end then, for fixed $\lambda_2$,  only one of the connected components in $\Cay(\Gamma)\setminus U_{\lambda_2}$ can be infinite. Up to action of $\Gamma\curvearrowright\Cay(\Gamma)$, there is just finite number of finite connected  components in $\Cay(\Gamma)$. Thus, we may take $L'$ such that $L'-L-10CDxN$ is greater than number of vertices of any finite connected component in any $\Cay(\Gamma)\setminus U_{\lambda_2}$.

Now, prove the desired for this $L'$ and for $\lambda_2$ being the beginning segment of $\lambda$ with $\length\lambda >L'$ satisfying conditions of our Lemma. Denote by $A$ the infinite connected component of $\Gamma\setminus U_{\lambda_2}$. For a Cayley graph vertex $v\in A\cap \Gamma$, let $\alpha(v)\in\lambda_2$ be such that $[v,\alpha(v)]$ is a geodesic perpendicular from $v$ to $\lambda_2$. 

By the choice of $L'$, geodesic segment $\lambda_1\setminus U_{\lambda_2}$ has length $\ge L'-L$ and thus cannot belong to a finite connected component of $\Cay(\Gamma)\setminus U_{\lambda_2}$, therefore, we may pick a  vertex  $v_1\in\lambda_1\cap A$. By Lemma \ref{lemma:perpendicular_stab2}, $\dist_{\Cay(\Gamma)}(\alpha(v_1),v_0)\le 2CDx$. Also, at $\lambda\setminus U$, there is a vertex $v_2$ with $\dist_{\Cay(\Gamma)}(v_2, \lambda_2)>L$ and we have $v_2\in A$ by the choice of $L'$ again. Since $\lambda$ is geodesic segment, we have $\alpha(v_2)=\ennd\lambda_2$.

Join $v_1$ and $v_2$ with a path in $A$. When $v\in A\cap \Gamma$ moves along this path by a distance $1$, that is, over an edge, then $\alpha(v)$ moves along $\lambda_2$ by a distance $\le 2CDx$. This is by Lemma \ref{lemma:perpendicular_stab}, otherwise we finish the proof of Lemma \ref{lemma:embed}. Since $\alpha(v)$ travels from a point near $\beg\lambda_2$ to $\ennd\lambda_2$, we arrive to the desired.
$\blacksquare$

\medskip

Now we may implement construction of geodesic segments $\lambda(\nu)$, $\nu$ ranges $\mathcal T$, in the order specified above but with (\ref{eq:fixed_step}) replaced by (\ref{eq:vary_step}). Let $\length\lambda(\nu)$ be depending only on level $k=0, 1, \dots, N-1$ of $\nu$ in $\mathcal T$, denote it by $L_k$. Pick $L_0, L_1, \dots, L_{N-1}$ such that $L_0\ge L_1'$, $L_k \ge L_{k+1}'$, $L_{k-1}\ge2L_{k+1}'$, $k=1,2,\dots, N-2$, $L_{N-1}\ge 4CDx$, and also such that $L_0\ge L_1\ge L_2\ge\dots\ge L_{N-1}$. (The latter inequality, in fact, follows from the construction of $L'$ in Lemma \ref{lemma:build_perpendicular} if we have $L_k\ge L_{k+1}'$.)

As above, we start with constructing $\lambda(\nu_0)$, where, recall, $\nu_0\in\mathcal T$ is the root. Let $\lambda(\nu_0)$ be half of a geodesic segment of length $2L_0\ge 2L_1'$. 
Put $\phi(\nu_0):=\beg\lambda(\nu_0)$. By Remark after  Lemma \ref{lemma:build_perpendicular}, we may build perpendiculars to $\lambda(v_0)$ with lengths $L_1$.  By Lemma \ref{lemma:build_perpendicular}, we may chose points $\phi(\nu_0):=\beg\lambda(\nu_0)$, $\phi(\nu_0r), \phi(\nu_0rl), \phi(\nu_0rl^{N-1})$ along $\lambda(\nu_0)$ 
and geodesic segments $\lambda(\nu_0r)$,  $\lambda(\nu_0rl)$,  $\lambda(\nu_0rl^{N-2})$ such that all the required conditions for these segments are satisfied.

Now, repeat this procedure for newly constructed segments of the form $\lambda(\nu)$. Let $k$ be the level of $\nu$ at the tree. We set $\phi(\nu) := \beg \lambda(\nu)$. We are going to apply Lemma~\ref{lemma:build_perpendicular} for $\lambda(\nu)$ to build perpendiculars of lengths $L_{k+1}$ with steps  in $[2CDx, 4CDx]$ along $\lambda(\nu)$. For the first assumption of Lemma~\ref{lemma:build_perpendicular}, it is enough that $L_k\ge L_{k+1}'$. 

We also need to check the second assumption in Lemma \ref{lemma:build_perpendicular}. To this end, notice that if $\mathcal T\ni \nu\neq \nu_0$ then there exists $\nu'\in\mathcal T$ such that $\nu=\nu'rl^j$ for some $j=0,1,2,\dots$. Segment $\lambda(\nu')$ is already built, and $\lambda(\nu)$ reversed is perpendicular to $\lambda(\nu')$. If $L_{k-1}>2L'_{k+1}$ then either of the two segments $[\phi(\nu'), \phi(\nu)]$ 
or  
$\lambda(\nu')\setminus[\phi(\nu'), \phi(\nu)]$ (both are subsets of $\lambda(\nu')$) can be taken as $\lambda_1$ in Lemma \ref{lemma:build_perpendicular}. We thus conclude that one may build geodesic perpendiculars along $\lambda(\nu)$ and define $\lambda(\nu r), \lambda(\nu rl), \dots, \lambda(\nu rl^{N-k-2})$ together with their beginnings $\phi(\nu r), \phi(\nu rl), \dots, \phi(\nu rl^{N-k-2})$ and also define leaf image $\phi(\nu rl^{N-k-1})$ such that (\ref{eq:vary_step}) will be held for  $\nu_1, \nu_2\in\{\nu, \nu r, \nu rl, \dots, \nu rl^{N-1-k}\}$ adjacent in the almost-binary tree.

To build leafs,  it is enough that  $L_{N-1}\ge 4CDx$.  We thus see that, 
under our choice of $L_k$, one may  successfully construct segments $\lambda(\nu)$ and vertices $\phi(\nu)$, all the required properties are held. If $\nu_1, \nu_2\in\mathcal T$ are adjacent then $\phi(\nu_1), \phi(\nu_2)$ belong to the same segment of the form $\lambda(\nu)$, it is seen from our construction. Let us also define $\phi$ at $\edge(\nu_1, \nu_2)$ of the almost-binary tree so as when $\nu'$ travels along the latter edge from $\nu_1$ to $\nu_2$ with unit speed then $\phi(\nu')$ travels with constant speed along $\lambda(\nu)$ from $\phi(v_1)$ to $\phi(v_2)$. We thus constructed a continuous $\phi$ from almost-binary tree to $\Cay(\Gamma)$, now also on edges of the former. Also, if $\nu\in\mathcal T$ is a leaf then we may define $\lambda(\nu):=[\phi(\nu), \phi(\nu)]$, a degenerate perpendicular; this will uniformize notation. 

Now we are going to catch a loop in the image of $\phi$ and then shrink it to have (\ref{eq:x_bilip}) and injectivity. Consider tuples $(\gamma_0, v)$ such that:
\begin{itemize}
	\item $\gamma_0$ is  curve in $\Cay(\Gamma)$ starting and ending at $\Gamma$, also $v\in\Gamma$.
	
	\item  $\gamma_0$ is injective except for the possibility that $\beg\gamma_0=\ennd\gamma_0$.
	
	\item $\dist_{\Cay(\Gamma)}(v, \beg\gamma_0)\le x$, $\dist_{\Cay(\Gamma)}(v, \ennd\gamma_0)\le x$.
	
	\item There exists a path $\gamma^{\mathcal{T}}$ in the almost-binary tree such that $\gamma_0=\phi(\gamma^{\mathcal T})$. Path $\gamma^{\mathcal{T}}$ here does not have to start or stop in $\mathcal T$ but may also have ends inside of edges of almost-binary tree.

\end{itemize}


Divide $\gamma_0$ by points of the form $\phi(\nu)$, $\nu\in\mathcal T$, belonging to $\gamma_0$. Let us call the closed non-degenerate arcs of such subdivision the \emph{sides} of $\gamma_0$. By the construction, all sides are geodesic segments. If two of  
sides $s_1, s_2$ are adjacent then one of them is, up to reverse of orientation, a perpendicular to another. This is true even if $s_1$ and $s_2$ are adjacent parts of some $\lambda(\nu)$. Let us call $s_1\cup s_2$ a \emph{corner} of $\gamma_0$ and also call \emph{corner point} the unique point in $s_1\cap s_2$. Denote by $s_b$ and $s_e$ the sides of $\gamma_0$ containing $\beg\gamma_0$ and $\ennd\gamma_0$, respectively. 

We impose one more  condition on tuple $(\gamma_0, v)$:

\begin{itemize}
	\item $\gamma_0$ has at least two  corners.
\end{itemize}


Let us call tuples $(\gamma_0, v)$ satisfying all the conditions above \emph{admissible}.

\begin{lemma}
	There exists at least one admissible tuple, provided that $$\card\mathcal B_{\Gamma}(1_\Gamma, 4CDxN) < 2^{N-1}.$$
\end{lemma}

\noindent {\bf Proof.} 
Consider injective paths $\gamma^{\mathcal T}$ in the almost-binary tree such that $\beg\gamma^{\mathcal T}\neq \ennd\gamma^{\mathcal T}$ but $\phi(\beg\gamma^{\mathcal T})=\phi(\ennd\gamma^{\mathcal T})$. They do exist since for a leaf $\nu\in\mathcal T$ image  $\phi(\nu)$ lies at $\mathcal B_{\Gamma}(\phi(\nu_0), 4CDxN)$ and there are $2^{N-1}$ images of leafs. 

Take the shortest $\gamma^{\mathcal T}$ as above. Then $\gamma_0:=\phi\circ {\gamma^{\mathcal T}}$ is injective except for the ends, for otherwise we may shorten $\gamma^{\mathcal T}$ by dropping a loop from it. This $\gamma_0$ has at least two corners because self-intersection is impossible at the sides of one corner. 
It remains to take $v:=\phi(\beg\gamma^{\mathcal T})=\phi(\ennd\gamma^{\mathcal T})$.
$\blacksquare$

\medskip

Now we conclude the proof of Lemma \ref{lemma:embed}. Among all the admissible tuples $(\gamma_0,v)$ take ones with minimal $\length\gamma_0$. Further, among the latter tuples, consider the one minimizing $\length([\beg\gamma_0,v]+[v,\ennd\gamma_0])$.

\begin{lemma}[on shrinking]
	\label{lemma:shrinking}
	For tuple $(\gamma_0,v)$ chosen as above, loop $$\gamma:=\gamma_0+[\ennd\gamma_0,v]+[v,\beg\gamma_0]$$ satisfies the requirements from Lemma \ref{lemma:embed}.
\end{lemma}

\noindent {\bf Proof.}
First we check condition (\ref{eq:x_bilip}). Assume that two points $p_1, p_2\in\gamma\cap \Gamma$ are such that $\dist_{\Cay(\Gamma)}(p_1, p_2)\le x$. Denote by $\dist_\gamma(p_1, p_2)$ the shortest distance along $\gamma$ between $p_1$ and $p_2$, here $\gamma$ is passed with unit speed.

If $p_1, p_2\in[\ennd\gamma_0,v]+[v,\beg\gamma_0]$ then 
$$\dist_\gamma(p_1, p_2)\le \length([\ennd\gamma_0,v]+[v,\beg\gamma_0])\le 2x\le Cx.$$

Check the case when $p_1, p_2\in\gamma_0$ but at least one of them is not in $\{\beg\gamma_0, \ennd\gamma_0\}$. If $p_1$ and $p_2$ lie on the same corner of $\gamma_0$ then, by Lemma \ref{lemma:perpendicular}, 
$\dist_{\gamma_0}(p_1, p_2)\le Cx$, the desired. If $p_1, p_2\in\gamma_0$ and do not lie at the same corner then tuple
$$
(\mbox{arc of }\gamma_0\mbox{ from }p_1\mbox{ to }p_2, p_1)
$$
is admissible and has a shorter curve, this contradicts our choice of tuple.

It remains to check the case 
\begin{align*}
&p_1\in\gamma_0\setminus\{\beg\gamma_0,\ennd\gamma_0\}, \\
&p_2\in([\ennd\gamma_0,v]+[v,\beg\gamma_0])\setminus\{\beg\gamma_0,\ennd\gamma_0\}.
\end{align*}
Without loss of generality, assume $p_2\in[v,\beg\gamma_0]$, $p_2\neq \beg\gamma_0$. Let $\gamma_0'$ be arc of $\gamma_0$ from $\beg\gamma_0$ to $p_1$. If 
$\gamma_0$ has at least two corner points inside of $\gamma_0'$ then tuple $(\gamma_0', p_2)$ is admissible and has a shorter curve which contradicts the choice of minimal tuple $(\gamma_0,v)$.

Let $s$ be side of $\gamma_0$ next to $s_b$, its first side. 
We need to check the case when $p_1\in s$ or $p_1\in s_b$. We have that one of $s_b$, $s$ is a perpendicular to another, up to orientation reverse. In both cases, arguing as in Lemma \ref{lemma:perpendicular}, we conclude that $\length\gamma_0'\le 6x$, $\dist_\gamma(p_1, p_2)\le 7x$. We then successfully check (\ref{eq:x_bilip}) provided that $C>7$.

Now, we also  have to check injectivity of loop $\gamma$. Curve $\gamma_0$ itself is injective except, possibly, for its ends. If $[\ennd\gamma_0,v]+[v,\beg\gamma_0]$ is not injective then we may construct an admissible tuple with the same $\gamma_0$ and smaller $\length([\ennd\gamma_0,v]+[v,\beg\gamma_0])$, a contradiction to the choice of tuple.

Assume, without loss of generality, that there is $p\in[v,\beg\gamma_0]\cap \gamma_0\cap \Gamma$, $p\neq\beg\gamma_0, \ennd\gamma_0$. 
If arc $\gamma_1$ of $\gamma_0$ between $\beg\gamma_0$ and $p$ passes at least two corner points of $\gamma_0$ then $(\gamma_1, p)$ is an admissible tuple with $\length \gamma_1<\length\gamma_0$, a contradiction. Similarly, if arc $\gamma_2$ of $\gamma_0$ from $p$ to $\ennd\gamma_0$ passes at least two corner points of $\gamma_0$ then $(\gamma_2, v)$ is an admissible tuple with $\length \gamma_2<\length\gamma_0$, a contradiction again. 

Thus, any of $\gamma_1$ and $\gamma_2$ has at most one corner point strictly inside it. (By the way, we have not excluded the possibility that $p$ is a corner point, and there are two more corner points at $\gamma_0$.) In this case, notice that 
$$
\dist_{\Cay(\Gamma)}(\beg\gamma_0, p)\le x, ~~  \dist_{\Cay(\Gamma)}(p,\ennd\gamma_0)\le 2x.
$$
Arguing as in Lemma \ref{lemma:perpendicular}, we conclude that $$\dist_{\gamma_0}(\beg \gamma_0, p)\le 3x, ~~ \dist_{\gamma_0}(p,\ennd \gamma_0)\le 6x.$$
But then $\length\gamma_0 \le 9x$ which is impossible if $C>9$ since $\gamma_0$ has at least a whole side of length $\ge 2CDx$, by the construction.
$\blacksquare$

\section{Spectral and homological argument}
\label{sec:spectra_homology}
In this Section, having already loop embedding given by Lemma \ref{lemma:embed}, we accomplish 

\medskip

\noindent {\bf Proof of Theorem \ref{th:exp_growth}}.
On $\ell_{0,c}^2(E)$, operator $\Delta_1$ is bounded and separated from zero. Thus, there exists $c<1$ such that, for any $n\in\mathbb N$, there exists polynomial $P_n$ of one variable with degree $n$ such that 
\begin{equation*}
	\left|P_n(\mathscr t)-1/{\mathscr t}\right|\le \const\cdot c^n, ~~ {\mathscr t}\in\spec(\Delta_1|_{\ell_{0,c}^2(E)});
\end{equation*}
here  $c$  depends on the size of spectral gap under consideration.
Thus, by Spectral Theorem,
\begin{equation*}
	\left\|P_n(\Delta_1)-\Delta_1^{-1}\right\|_{\ell_{0,c}^2(E)\to \ell_{0,c}^2(E)}\le \const\cdot c^n.
\end{equation*}

Laplacian $\Delta_1$ is a local operator, namely, $\supp \Delta_1f$ is contained in $D$-neighborhood of $\supp f$  (in $\Cay(\Gamma)$-metric) for any $f\in\ell^2_{0,c}(E)$. Indeed, for some $e\in E$, cochain $\Delta_1 \mathds 1_{e}$ is supported by edges which belong to  faces from $\Cay^{(2)}(\Gamma)$ containing $e$. We derive that $P_n(\Delta_1)f$ is supported by $(D\cdot n)$-neighborhood of $\supp f$. 

Let $\gamma$ be a loop of length $\mathscr L$  provided by Lemma \ref{lemma:embed} for $x=2Dn+2D+1$, $C:=21$, and, further,  $f_\gamma\in\ell_{0,c}^2(E)$ be closed $1$-cochain (or, rather, chain) given by $\gamma$. 
Then, $\supp f_\gamma=\gamma(\mathbb T_{\mathscr{L}})$. Recall that on $\ell_{0,c}^2(E)$ we have $\Delta_1=\partial d$. Consider $1$-cochain
\begin{equation*}
	g:=f_\gamma - \Delta_1 P_n(\Delta_1)f_\gamma = f_\gamma - \partial d  P_n(\Delta_1)f_\gamma = \Delta_1\left(\Delta_1^{-1}-P_n(\Delta_1)\right)f_{\gamma}.
\end{equation*}

By the choice of $P_n$, we have 
\begin{equation}
	\label{eq:g_estim_l2}
	\|g\|_{\ell_{0,c}^2(E)} \le \const\cdot c^n \cdot \|f_\gamma\|_{\ell^2_{0,c}(E)}=\const\cdot c^n \cdot\sqrt{\length\gamma}.
\end{equation}
Let $U$ be $(Dn+D)$-neighborhood of $\gamma$. We have that  $f_\gamma-g$ bounds in $U$, 
that is, is $\partial$ of a $2$-cochain supported by $U$. 


Thus, we can spread $f_\gamma$ in $U$ such that the resulted $1$-cochain $g$ is homologic to the original $f_\gamma$ but this new $g$ is exponentially small with respect to $f_\gamma$ in $\ell^2$-norm. Let us take some informal consideration. If we could somehow speak about "sections $\sigma$ of $U$ perpendicular to $\gamma$"{} then, for each such $\sigma$ "flows" of $f_\gamma$ and of $g$ through $\sigma$ should coincide due to homology between $1$-(co)chains. The flow of $f_\gamma$ through $\sigma$ is $1$. By summing up over all $\sigma$, this would allow to estimate $\|g\|_{\ell^1(E)}\gtrsim \length\gamma$. Together with $\ell^2$-smallness of $g$, this leads to lower estimate for $\card\supp g\le \card U$ which is enough for us. 

We return to the formal argument. On \emph{vertices} passed by $\gamma$, introduce cyclic coordinate 
$$
\varphi \colon \Gamma\cap \gamma\to \mathbb Z\mmod \length\gamma.
$$
We are going to extend this coordinate to $U$ to obtain a uniformly locally Lipschitz multivalued function. 
In our construction, we make use of metric conditions on the curve. 

Pick any $v\in U\cap \Gamma$. Drop a perpendicular from $v$ to $\gamma$, 
namely, let $\alpha(v)\in\gamma(\mathbb T_{\mathscr L})\cap \Gamma$ be any of vertices passed by  $\gamma$ closest to $v$ in graph metric. Put $\varphi(v):=\varphi(\alpha(v))\in \mathbb Z\mmod \length\gamma$. 

Let $v_1, v_2\in U\cap \Gamma$ be vertices  adjacent in $\Cay(\Gamma).$ Then $$\dist_{\Cay(\Gamma)}(\alpha(v_1), \alpha(v_2)) \le 1+2Dn+2D$$ which is $x$. Then, by the choice of $\gamma$,
\begin{equation}
	\label{eq:phi_diff}
	\dist_{\mathbb Z \mmod \length\gamma}(\varphi(v_1), \varphi(v_2)) \le Cx.
\end{equation}
We define $\psi\colon E\to \R$, roughly speaking, to be $d\varphi$ strictly inside of $U\cap E$. More carefully, if $e\in E$ and either $\beg e\notin U$ or $\ennd e\notin U$  then define $\psi(e)$ arbitrarily. Otherwise, let $\psi(e)$ be a real number which belongs to $\left(\varphi(\alpha(\ennd e))-\varphi(\alpha(\beg e))\right)+ \mathbb{Z}\cdot\length \gamma$ and has the least absolute value over this set. We, in particular, have estimate $|\psi(e)|\le Cx$ if $\beg e, \ennd e\in U$.

We claim that $\langle \psi, \partial d P_n(\Delta_1)f_\gamma\rangle_{\ell^2(E)}=0$. Indeed, since $\supp P_n(\Delta_1)f_\gamma$ is finite, the latter is $\langle d\psi, dP_n(\Delta_1)f_\gamma\rangle_{\ell^2(F)}$. Let $\sigma \in \supp dP_n(\Delta_1)f_\gamma\subset F$ be a face. Since $\supp P_n(\Delta_1)f_\gamma$ lies in $(Dn)$-neighborhood of $\gamma$, all edges in $\sigma$ belong to $U$. When $v$ ranges vertices from $\sigma$, 
$\alpha(v)$ belongs to an arc in $\gamma(\mathbb T_{\mathscr L})$ having length $\le D\cdot Cx$. This is because $\sigma$ has $\le D$ vertices and by (\ref{eq:phi_diff}) for adjacent vertices. Since $\length \gamma>2CDx$, we conclude that $d\psi(\sigma)=0$. This implies that $\langle \psi, \partial d P_n(\Delta_1)f_\gamma\rangle_{\ell^2(E)}=0$.


If some $e$ is passed by curve $\gamma$ then $f_\gamma(e)=\psi(e)\in\{-1,+1\}$ due to injectivity of~$\gamma$. Thus we notice that $\langle f_\gamma,\psi\rangle_{\ell^2(E)}=\length\gamma$. 

Now write 
\begin{multline*}
	\length\gamma=\langle f_\gamma, \psi\rangle_{\ell^2(E)}=\langle g+\partial d P_n(\Delta_1)f_\gamma, \psi\rangle_{\ell^2(E)}=\langle g, \psi\rangle_{\ell^2(E)} \le Cx \cdot \|g\|_{\ell^1(U\cap E)}\le \\ \le Cx \cdot \sqrt{\card\left(U\cap E\right)}\cdot \|g\|_{\ell^2(U\cap E)} \le 
	Cx\cdot\sqrt{\card\left(U\cap E\right)}\cdot c^n \cdot \|f_\gamma\|_{\ell^2(E)} = \\=
	Cx\cdot\sqrt{\card\left(U\cap E\right)}\cdot c^n \cdot \sqrt{\length\gamma}.
\end{multline*}
By (\ref{eq:g_estim_l2}), we then have 
$$
\card\left(U\cap E\right) \ge \const\cdot\left(\frac{1}{c^2}\right)^n\cdot \length\gamma\cdot \frac{1}{n^2}.
$$
This immediately implies exponential growth of balls in $\Gamma$ since $$\card(U\cap E)\le \const\cdot\length\gamma\cdot \card(\mathcal B_\Gamma(1_\Gamma,n)). ~ \blacksquare$$

\medskip

\noindent {\bf Remark.} Boundedness of the resolvent $(\Delta_1|_{\ell^2_{0,c}(E)})^{-1}$ can be reformulated as follows: \emph{there exists $\mathscr C<+\infty$ such that for any $f\in\ell_{0,c}^2(E)$ there exists 
$h\colon F\to\R$ such that 
\begin{equation}
	\label{eq:ell2_film}
\partial h=f \mbox{ and }\|h\|_{\ell^2(F)}\le \mathscr C\cdot \|f\|_{\ell^2(E)}.
\end{equation}} Indeed, in the case of presence of a spectral gap, one takes $h:=d\Delta_1^{-1}f$. To prove the opposite, we notice that the assumption on existence of $h$ as above implies that 
$$
\partial\colon \ell^2(F)\ominus_{\ell^2(F)} \Ker(\partial|_{\ell^2(F)}) \to \ell^2_{0,c}(E)
$$
is an open operator and thus a bijection, this implies that $\Delta_1=\partial\partial^*$ is also bijective on $\ell^2_{0,c}(E)$.  Notice also that $2$-cochains from $\ell^2(F)\ominus_{\ell^2(F)}\Ker(\partial|_{\ell^2(F)})$ minimize $\ell^2$-norm with prescribed $\partial$.

If $U$ is as in the argument above then, for $h$ as in (\ref{eq:ell2_film}), we have that $h$ has to have a significant $\ell^1$-mass on boundary of $U$, this also means that either $\|h\|_{\ell^2(F)}$ or boundaries of $U$ are large.

Consider Euclidean space $\R^3$. Though Laplacian here is not invertible, one may consider analog of $d(\Delta_1)^{-1}\colon\ell_{0,c}^2(E)\to \ell^2(F)$. Let $f\colon \R^3\to\R^3$ be a $C_0^\infty$-vector field. Denote by $\star$ the distributional convolution. Then magnetic \emph{Biot--Savard field} $\BS^f=\curl(\frac{1}{4\pi|x|}\star f)$ solves $\curl h=f$. Approximation properties of such potentials were studied in \cite{HP96}, \cite{MH98} which partially motivated the current work.  


\medskip

\noindent {\bf Remark.} One may expect application of \emph{frame operator} instead of $\Delta_1^{-1}$ and its polynomial approximation. But frame techniques itself is still out of use here.

\medskip


\medskip

\section{Examples}

\label{sec:examples}

In this Section, we concern Riesz systems.  The following definition agrees to the one given at Section \ref{sec:maueutic}:

\begin{define}
	Let $H$ be a separable Hilbert space,  $f_1, f_2, \dots\in H$. We say that $\{f_n\}_{n\in\mathbb N}$ is a \emph{Riesz system} if for any \emph{finitely supported} sequence of coefficients $\{a_n\}_{n\in\mathbb N}\subset \R$, we have the following \emph{almost-orthogonality} relation:
	\begin{equation}
		\label{eq:almost_orth}
		C^{-1}\cdot \sum\limits_{n\in\mathbb N} a_n^2 \le 
		\left\|\sum\limits_{n\in\mathbb N} a_n f_n\right\|^2 \le 
		C\cdot \sum\limits_{n\in\mathbb N} a_n^2 
	\end{equation}
	with some $C\in(0,+\infty)$ not depending on $\{a_n\}_{n\in\mathbb N}$.
	%
	
	A \emph{Riesz basis} in $H$ is a complete Riesz system therein.
\end{define}

Thus, Riesz system is a Riesz basis in its norm-closed linear span.

As we indicated in Section \ref{sec:maueutic}, existence of a Riesz basis in $\ell^2_{0,c}(E)$ with diameters bounded from the above implies that $1$-Laplacian has a spectral gap on this space.

\subsection{$\mathbb Z^m$}

By now, we already know that there is no spectral gap for $\Delta_1$ at $\mathbb Z^m$, $m\ge 2$. This follows from Theorem \ref{th:exp_growth}, but can be seen immediately. Consider, for simplicity, $\mathbb Z^2$ with standard generators. For $n\in\mathbb N$, let $Q_n\subset \R^2$ be square with center in $(0,0)$ and side $2n$. Then, for any $N=1,2,\dots,$ put $f:=\sum_{n=1}^N \partial Q_n$. This $f$ is indeed a closed cochain on edges of $\mathbb Z^2\subset \R^2$ but, in $\ell^2$-norm, its curl is negligible with respect to $f$ itself for $N$ large.  

It follows that there does not exits a nice Riesz basis in $\ell^2_{0,c}(E)$. But we may be curious on Riesz systems. 

\begin{predl}
	Consider $\mathbb Z^m$ equipped with generators $\mathscr e_j=(0, \dots, 0, 1, 0, \dots, 0)$ with $1$ at $j$th position, $j=1, \dots, m$. 
	
	In the corresponding Cayley graph, $\mathbb Z^m$-shifts of any finitely supported closed $f\in \ell^2_{0,c}(E)$ do not form a Riesz system.
\end{predl}

Indeed, for any $j=1, \dots, m$, we have $\sum_{z\in\mathbb Z^m} f(\edge(z, z+\mathscr e_j))=0$ due to solenoidality of $f$. Therefore, for $R>0$ large enough, Riesz system condition fails for 
$$
\sum\limits_{z\in[-R,R]^m\cap\mathbb Z^m}f(\cdot -z)
$$
since the latter sum does not vanish only near the boundary of $[-R,R]^m$ and cancels at the rest of the lattice. 

In fact, we may say a bit more: 

\begin{predl}
	\label{predl:Zd_noinvariant}
	For any $m=2,3,\dots$, there is no translation invariant Riesz system in $\ell^2_{}(E)$ consisting of closed $1$-cochains belonging to $\ell^1(E)$. 
\end{predl}

Proof is obtained by a straightforward application of Fourier analysis and is given at the Appendix.

%
%

If we duplicate at least one generator, say add $\mathscr e_1':=\mathscr e_1$ to generating system then the conclusion of the Proposition above will not be true. Indeed, consider $f$ given by cycle $\edge(0,\mathscr e_1)$,  $\edge(\mathscr e_1,\mathscr e_1+\mathscr e_2)$, $\edge(\mathscr e_1+\mathscr e_2,\mathscr e_1+\mathscr e_2-\mathscr e_1')$, $\edge(\mathscr e_2,0)$ with the third edge given by $\mathscr e_1'$ reversed. Then $\mathbb Z^m$-shifts of $f$ do not cancel and form a Riesz system.




\subsection{Fundamental group of 2D sphere with $\mathscr g\ge 2$ handles} 

\label{subsec:Laplace}

Consider $2$-dimensional sphere with $\mathscr g\ge 2$ handles. Pick a canonical system of generators for its fundamental group $\Gamma$, that is, we write $\Gamma$ via presentation
\begin{equation*}
\Gamma = \langle a_1, b_1, \dots, a_{\mathscr g}, b_{\mathscr g} \mid a_1 b_1 a_1^{-1} b_1^{-1} \dots a_{\mathscr g} b_{\mathscr g} a_{\mathscr g}^{-1} b_{\mathscr g}^{-1}=1_{\Gamma}\rangle.
\end{equation*}

It is well-known that $\Gamma$ is quasiisometric to hyperbolic plane $\HH$. So, consider canonical tessellation $T$ of  $\HH$ corresponding to $\Gamma$.  
Any element of $T$ is a hyperbolic polygon in $\HH$ with $4{\mathscr g}$ vertices and is a fundamental domain of an action $\Gamma\curvearrowright\HH$, so 
one may understand $\Gamma$ as a subgroup in the group of orientation-preserving isometries of $\HH$. 
%
In terms of tessellation, vertices of $\Cay(\Gamma)$ 
are identified to polygons in $T$; by the choice of generators system,  two such elements are adjacent in  
$\Cay(\Gamma)$ iff they have a common edge as polygons in $\HH$.

Let $\mathscr V$ be the set of vertices of polygons from $T$. Pick any $v\in\mathscr  V$. 
Draw loop in $\widetilde{\Cay}(\Gamma)$ 
whose vertices are $\{F\in T\colon v \mbox{ is a vertex of } F\}$ and are ordered, say, in the cyclic positive direction as seen from $v$. By writing $\pm 1$ on edges of this cycle and $0$ on all the other edges of $\Cay(\Gamma)$, we arrive to cycle $\gamma_v\in \ell_{0,c}^2(E)$.

\begin{predl}
Family $\{\gamma_v\}_{v \in \mathscr V}$ is a Riesz basis in $\ell^2_{0,c}(E)$.

Thus, 
$\Gamma$ endowed with any other system of generators  has a coexact spectral gap for $1$-Laplacian. 
\end{predl}

{\noindent \bf Remark.} For author, it is still an open question whether existence of a Riesz basis survives under change of generators even in this case. 

\medskip

{\noindent \bf Remark.} The proof below is based on a discrete version of Hodge star, so is specific for dimension $2$ of the hyperbolic plane.

\medskip

\noindent {\bf Proof.}
Clearly,  $\{\gamma_v\}_{v \in \mathscr V}$ is complete in $\ell_{0,c}^2(E)$. 
We need to check that $\{\gamma_v\}_{v \in \mathscr V}$ is a Riesz system in $\ell^2(E)$. Upper estimate in Riesz condition (\ref{eq:almost_orth}) is simple because all cycles $\gamma_v$ have number of edges bounded from the above (in fact, the same). 

To prove the lower estimate pick finitely supported family of coefficients $\{\mathscr a_v\}_{v\in \mathscr V
}$. Edges of $\Cay(\Gamma)$ are identified with those of its dual, that is, to edges of polygons from~$T$. We thus need to check that 
\begin{equation*}
\sum\limits_{e \mbox{\scriptsize{ -- an edge of tessellation }}T}\left(\mathscr a_{\mathrm{end}(e)}-\mathscr a_{\mathrm{begin}(e)}\right)^2 \ge \const \cdot \sum\limits_{v\in \mathscr V} \mathscr a_v^2.
\end{equation*}
But this, by Kesten Theorem, follows from non-amenability of $\Gamma$. $\blacksquare$

\medskip

Another possible way to solve the question on spectral gap on $\Gamma$ is to reduce it to the similar question on hyperbolic plane $\HH$. In the case of manifold, we have to check the estimate $\|\omega\|_{L^2}\lesssim \|d\omega\|_{L^2}$ for compactly supported $1$-form $\omega$ with zero codifferential. To prove the equivalence of existence of spectral gap in discrete and in continuous settings, a de Rham-type Theorem, one has to repeat proof from \cite{M08} which goes back to A. Weil (double complex and Whitney formula are involved). At $\HH$, one checks the existence of a coexact spectral gap either directly or using spectral decomposition from \cite{Do81}.

\section{Appendix: some technical proofs}

\label{section:appendix}

\noindent {\bf Proof of Proposition \ref{predl:general_invariance}.} \emph{Spectral gap  at $G_1$ implies spectral gap at  $G_2$.} First suppose that $G_1$, the smaller graph, possesses spectral gap; let $D>0$ be as in the construction of $G_1^{(2)}$. We are going to check the spectral gap property for  $G_2^{(2)}$ built  with  
$$
D' := \max\left(D, \sup\limits_{e\in E_2\setminus E_1} \dist_{G_1}(\beg e, \ennd e)+1\right).
$$
Pick $f\colon E_2\to\R$ 
compactly supported with $\partial_{G_2}f=0$. Put $f_1 := f|_{E_1}$. Then 
\begin{equation*}
	\|df_1\|_{\ell^2(F^{}G_1)} \le \|df\|_{\ell^2(F^{}G_2)}.
\end{equation*}
(Recall that $FG_{1,2}$ are the sets of faces in $2$-dimensional complexes obtained from the corresponding graphs.) 

%
We are not permitted to apply (\ref{eq:rot_estim}) for $G_1$ and $f_1$ because $\partial_{G_1}f_1$ generally does not vanish. Thus define $f_2 \in \ell^2_{0,c}(E_1)$ as 
$f_2 := \pr_{\ell^2_{0,c}(E_1)} f_1$ where $\pr$ denotes the orthogonal projection.
By 
Hodge decomposition, we may write 
\begin{equation}
	\label{eq:f_2_def}
	f_2 = f_1 - d_{G_1}w - g
\end{equation}
for some $w\colon V\to \R$ with $\Delta_{0,G_1}w=0$ and $g\in\clos_{\ell^2(E_1)}{\{d_{G_1}u\colon u\in \ell_0(V)\}}$ where $\ell_0(V)$ is the space of finitely supported vertex functions. Assumptions of our Proposition imply that $g=g_1|_{E_1}$ for some $g_1\colon E_2\to\R$ belonging to $\clos_{\ell^2(E_2)}{\{d_{G_2}u\colon u\in \ell_0(V)\}}$. Therefore, if we put 
$$
f_3 := f-d_{G_2}w-g_1\colon E_2\to\R
$$
then $f_2=f_3|_{E_1}$. Notice that $f\bot_{\ell^2(E_2)} d_{G_2}w+g_1$  since $f\in \ell^2_{0,c}(E_2)$. Then, under our choice of $D'$,
\begin{equation*}
	\|f\|_{\ell^2(E_2)} \le \|f_3\|_{\ell^2(E_2)}\le\const\cdot\left(\|f_3\|_{\ell^2(E_1)}+\|d_{G_2}f_3\|_{\ell^2(F^{}G_2)}\right).
\end{equation*}
The latter is because any edge from $E_2\setminus E_1$ is a part of a loop of length $\le D'$ in $G_2$ with all the other edges in $E_1$.

We have $d_{G_1} f_2 = d_{G_1}f_1$ by (\ref{eq:f_2_def}) and continuity of $d_{G_1}\colon \ell^2(E_1)\to \ell^2(F^{}G_1)$. Similarly, $d_{G_2}f_3=d_{G_2}f$. Now we are finally ready to make use of (\ref{eq:rot_estim}) for $G_1$ and $f_2$ to write  
\begin{multline*}
	\|f\|_{\ell^2(E_2)} 
	\le 
	\const\cdot\left(\|f_3\|_{\ell^2(E_1)}+\|d_{G_2}f_3\|_{\ell^2(F^{}G_2)}\right)
	=\\=
	\const\cdot\left(\|f_2\|_{\ell^2(E_1)}+\|d_{G_2}f\|_{\ell^2(F^{}G_2)}\right)
	\le
	\const\cdot\left(\|d_{G_1} f_2\|_{\ell^2(F^{}G_1)}+\|d_{G_2}f\|_{\ell^2(F^{}G_2)}\right)
	=\\=
	\const\cdot\left(\|d_{G_1} f_1\|_{\ell^2(F^{}G_1)}+\|d_{G_2}f\|_{\ell^2(F^{}G_2)}\right) \le \const\cdot \|d_{G_2}f\|_{\ell^2(F^{}G_2)}.
\end{multline*}
This proves the existence of a spectral gap at  $G_2$.

\medskip

\medskip

\noindent \emph{Spectral gap for $G_2$ implies spectral gap  for $G_1$.} 
Here, the argument is generally similar. Let $D'>0$ be the constant implied in the construction of $G_2^{(2)}$.  Now we put
$$
D := D'\cdot\left(1+\sup\limits_{e\in E_2\setminus E_1} \dist_{G_1}(\beg e, \ennd e)\right);
$$
we are going to check that if we build $G_1^{(2)}$  with this $D$ then we will have a spectral gap. 

Pick $f\in \ell_{0,c}^2(E_1)$, that is, with $\partial_{G_1}f=0$. For any $e\in E_2\setminus E_1$, pick any simple curve $\gamma_e$ in $\tilde G_1$ starting in $\beg e$ and ending at $\ennd e$ and of length $\dist_{G_1}(\beg e, \ennd e)$. Define $f_1\colon E_2\to \R$ such that $f_1|_{E_1}=f$ and, for $e\in E_2\setminus E_1$, $f_1(e)$ equals $\sum\limits_{e'\in \gamma_e} \pm f(e')$; here we take "$+$" sign if $\gamma_e$ passes $e'$ in its direction in $G_1$, and we take "$-$" sign otherwise. For such $f_1$ and $D$ as defined we see  that $\|d_{G_2}f_1\|_{\ell^2(F^{}G_2)} \le \const\cdot\|d_{G_1}f\|_{\ell^2(F^{}G_1)}$.

Write orthogonal decomposition 
\begin{multline*}
	f_1\in \ell^2(E_2) = \ell^2_{0,c}(E_2)\oplus_{\ell^2(E_2)} \clos_{\ell^2(E_2)}{\{d_{G_2}v\colon v\in \ell_0(V)\}}\oplus_{\ell^2(E_2)} 
	\\ 
	\oplus_{\ell^2(E_2)}\{d_{G_2}w\mid w\colon V\to \R,\, \Delta_{0,G_2}w=0,\,d_{G_2}w\in\ell^2(E_2)\}.
\end{multline*}
Let $f_2$ be the projection of $f_1$ to $\ell^2_{0,c}(E_2)$. We have $d_{G_2}f_2=d_{G_2}f_1$. Further, $f_1-f_2\bot_{\ell^2(E_2)}\ell^2_{0,c}(E_2)$, thus,  $f_1-f_2=d_{G_2}u$ for some $u\colon V\to\R$ with $d_{G_2}u\in\ell^2(E_2)$.
Then, $(f-f_2)|_{E_1}=(f_1-f_2)|_{E_1}=d_{G_1}u$. Since $f\in \ell^2_{0,c}(E_1)$, we have $f\bot_{\ell^2(E_1)} d_{G_1}u$. Thus, $\|f_2\|_{\ell^2(E_1)} \ge \|f\|_{\ell^2(E_1)}$. 

Now, using (\ref{eq:rot_estim}) for $G_2$ and $f_2\in\ell_{0,c}^2(E_2)$, we may write
\begin{multline*}
	\|f\|_{\ell^2(E_1)}\le\|f_2\|_{\ell^2(E_1)}\le
	\|f_2\|_{\ell^2(E_2)}\le \const\cdot\|d_{G_2}f_2\|_{\ell^2(F^{}G_2)}
	=\\=
	\const\cdot\|d_{G_2}f_1\|_{\ell^2(F^{}G_2)}\le 
	\const\cdot\|d_{G_1}f\|_{\ell^2(F^{}G_1)},
\end{multline*}
the desired. Proof is complete. $\blacksquare$

\medskip

\medskip

\noindent {\bf Proof of Proposition \ref{predl:scf_def_equivalence}.} Let 
$\mathcal L$ be the set of all simple loops in $\tilde G$ with lengths $\le D$. Any $\gamma \in \mathcal L$ defines a cycle $f_\gamma \in \ell^2_{0,c}(E)$. Any $g_n$, $n\in\NN$, can be $\ell^1$-convexly decomposed into the latter cycles:
\begin{equation*}
	g_n = \sum\limits_{\gamma \in \mathcal L} a_{n,\gamma} f_\gamma 
\end{equation*}
with some $a_{n, \gamma}\in\R$ such that 
\begin{equation}
	\label{eq:decomposition_convexity}
	\|g_n\|_{\ell^1(E)} = \sum\limits_{\gamma\in\mathcal L} |a_{n,\gamma}|\cdot \length\gamma.
\end{equation}
Also, if $a_{n,\gamma}\neq 0$ then $\supp\gamma\subset \supp g_n$. This is an elementary version of S.K. Smirnov Decomposition Theorem. 
We have that 
\begin{equation}
	\label{eq:decomposition_local_finiteness}
	C:=\sup\limits_{n\in\NN}\card\{\gamma\in\mathcal L\colon a_{n,\gamma}\neq 0\} <+\infty.
\end{equation}

We are going to show that $\{f_\gamma\}_{\gamma\in\mathcal L}$ is a frame in $\ell^2_{0,c}(E)$. The upper estimate from the frame definition is held automatically. The goal is to establish the lower one. Pick $y\in \ell^2_{0,c}(E)$. By Cauchy--Bunyakowskiy--Schwartz inequality, we have
\begin{equation*}
	\sum\limits_{n\in\NN} \langle g_n, y\rangle^2 =
	\sum\limits_{n\in\NN} \left\langle \sum\limits_{\gamma \in \mathcal L} a_{n,\gamma} f_\gamma, y\right\rangle^2
	\le 
	C\cdot \sum\limits_{\gamma\in\mathcal L}\left(\sum\limits_{n\in\NN}a_{n,\gamma}^2\right)\cdot \langle f_\gamma, y\rangle^2.
\end{equation*}
Here, $C<+\infty$ is the constant from (\ref{eq:decomposition_local_finiteness}).
Since we know the lower frame estimate for $\{g_n\}_{n\in\NN}$, it is enough to show that $\sup\limits_{\gamma\in\mathcal L}\sum\limits_{n\in\NN}a_{n,\gamma}^2<+\infty$. 

By uniform boundedness of supports of $g_n$ and by (\ref{eq:decomposition_convexity}), 
$
|a_{n,\gamma}| \le \const \cdot \|g_n\|_{\ell^2(E)}
$ with some constant not depending on $n$ and $\gamma$. In cycle decomposition, one has that  $\supp\gamma\subset \supp f_n$ if $a_{n,\gamma}\neq 0$;  then it is enough to prove that, for $e$ ranging $E$, the sum
\begin{equation*}
	\sum\limits_{n\colon e\in\supp g_n}\|g_n\|^2_{\ell^2(E)}
\end{equation*}
is bounded from the above uniformly by $E$.

To this end, we make use of the upper frame estimate for $\{g_n\}_{n\in\NN}$. 
Due to boundedness of degrees in $G$, there are, up to an isometry, just finite number of configurations of $\mathcal B_G(e,D)$ when $e$ ranges $E$. In any of such configurations, write upper frame estimate 
$$
\sum\limits_{n\colon \supp g_n \subset \mathcal B_G(e,D)} \langle y,g_n\rangle^2\le \const \cdot \|y\|^2_{\ell^2(E)}
$$
and average it over unit sphere in the space of closed $1$-cochains $y$ supported by $\mathcal B_G(e,D)$. This leads to the desired. Proof is complete.
$\blacksquare$

\medskip

\medskip

\noindent {\bf Proof of Proposition \ref{predl:Zd_noinvariant}.} It is enough to show the following: \emph{for any  $f\in \ell^1(E)$ with $\partial f=0$
	, the set of shifted cycles $\{f(\cdot -u)\}_{u\in \mathbb Z^m}$ is not a Riesz system in $\ell^2(E)$}. 

The upper estimate from Riesz system condition is obvious for $f\in \ell^1(E)$ just by Young convolution inequality for $\ell^1(\mathbb Z^m)* \ell^2(\mathbb Z^m)$. We are going to disprove the lower Riesz system estimate.

Write $f$ in coordinates: for $v\in \mathbb Z^m$ and $j=1,2,\dots, m$, put $$f_j(v) := f(\edge(v, v+\mathscr e_j)).$$
Pass to the dual group. If, in coordinate notation, $z=(z_1, \dots, z_m)\in\mathbb T^m$, $v=(v_1, \dots, v_m)\in\mathbb Z^m$ then write $z^v := z_1^{v_1}\cdot\dots\cdot z_m^{v_m}$.  For $z \in \mathbb T^m$, put $g_j(z) := \sum\limits_{v\in\mathbb Z^m} f_j(v) z^v$.

$\ell^2(E)$ is the same as $\left(\ell^2(\mathbb Z^m)\right)^m$. Riesz system condition for the family of all $\mathbb Z^m$-shifts of $f$ is taken by Fourier transform to the following: \emph{the set of all $m$-vector functions $(z^v g_1(z), \dots, z^v g_m(z))\in \left(L^2(\mathbb T^m)\right)^m$ with $v$ ranging $\mathbb Z^m$ is a Riesz system in $\left(L^2(\mathbb T^m)\right)^m$}. Write Riesz system condition for the latter system and for some family of coefficients $\{a(v)\}_{v\in\mathbb Z^m}\subset\mathbb C$. (The original system of shifts of $f$ was real-valued. It is no matter whether to consider real or complex coefficients in its Riesz condition. The same thus is true after Fourier transform.) We see that if $S(z) = \sum\limits_{v\in\mathbb Z^m} a_v z^v$ ($z\in\mathbb T^m$) then $\max\limits_{j=1, \dots, m} \|S(z) g_j(z)\|_{L^2(\mathbb T^m)}\asymp\|S\|_{L^2(\mathbb T^m)}$, the two-sided inequality with constants not depending on $S$; the latter holds for any $S\in L^2(\mathbb T^m)$. This is possible if and only if $|g_1|+\dots +|g_m|$ is bounded from the above and separated from zero on $\mathbb T^m$. We are going to disprove the lower part of the latter condition.

$f$ is a cycle. Hence we have 
\begin{equation}
	\label{eq:Fourier_div}
	\sum\limits_{j=1}^m g_j(z) (1-z_j)=0.
\end{equation}
Also, $f\in\ell^1(E)$, then each $g_j\in C(\mathbb T^m)$. In (\ref{eq:Fourier_div}) put $z=(z_1, \dots, z_m)$ 
with $z_j\neq 1$, $z_{j'}=1$ for $j'\neq j$. We see that $g_j(z)=0$ for such $z$. Hence $g_j(1, \dots, 1)=0$ by continuity. Again by continuity, we conclude that $|g_1|+\dots +|g_m|$ cannot be separated from zero near $z=(1, \dots, 1)$. Proof is complete.
$\blacksquare$

{\small
}

\end{document}